\documentclass{amsart}
\pdfoutput=1

\usepackage{amsmath}
\usepackage{amssymb}
\usepackage{graphicx}
\usepackage{url}

\newtheorem{theorem}{Theorem}[section]

\theoremstyle{definition}

\theoremstyle{remark}

\numberwithin{equation}{section}

\newfont{\hiera}{cmsy10 scaled 2488}
\newfont{\hierb}{cmsy10 scaled 1728}
\newfont{\hierc}{cmsy10 scaled 1200}

\newcommand{\Bigast}{
\mathop{\vphantom{\sum}\lower2.5pt\hbox{\hiera\char3}}}%

\newcommand{\Bigtimes}{
\mathop{\vphantom{\sum}\lower2.5pt\hbox{\hiera\char2}}}%

\begin{document}

\title{Eleven Unit Distance Embeddings of the Heawood Graph}

\author{Eberhard H.-A. Gerbracht}
\curraddr{Bismarckstra\ss e 20, D-38518 Gifhorn, Germany}
\email{e.gerbracht@web.de}

%    General info

\date{}


\begin{abstract}
In this note we present eleven unit distance embeddings of the Heawood graph, i.e.\ the point-line incidence graph of the finite projective plane of order two, by way of pictures and 15 digit approximations of the coordinates of the vertices. These together with the defining algebraic equations suffice to calculate arbitrarily exact approximations for the exact embeddings, and so to show that the Heawood graph indeed is a unit-distance graph. Thus we refute a ``suspicion'' of V.\ Chvatal from 1972.  
\end{abstract}

\maketitle

\section{Introduction}
In an informal joint collection of ``selected combinatorial research problems'' V.\ Chvatal in 1972 asked for characterizations of unit-distance graphs\footnote{See problem 21 in \cite{CombinatorialProblems}.}. These are graphs which can be embedded into the Euclidean plane in such a way that vertices correspond to points in the plane and adjacent vertices are exactly at distance one from each other, i.e., points corresponding to adjacent vertices can be connected to each other by straight unit length line segments\footnote{Synonymously we speak of unit distance embeddable graphs.}.

To make some headway Chvatal especially asked if the point-line incidence graphs of finite projective planes are unit distance embeddable -- and voiced the suspicion that they were {\sl not}. Even for the smallest of these,  associated to the projective plane of order two -- the so-called {\sl Heawood \hbox{Graph --}}  which is a graph with $2\cdot (2^2+2+1) = 14$ vertices and $(2+1)\cdot (2^2+2+1) =21$ edges, this question has remained unanswered until today. 

Partially inspired by this author's analysis of the Harborth graph \cite{GerbrachtHarborth}, which used dynamic geometry software and computer algebra to give final proof of the unit distance embeddability of the Harborth graph, in \cite{HarrisUDE} M.\ Harris  described a general strategy to approach this problem in case of the Heawood graph, leaving out all necessary calculations. Even though he thus sketched the basic idea of how to find a unit distance embedding for the Heawood graph, he did not provide any example.

In this note we present eleven unit distance embeddings of the Heawood graph by way of pictures and 15 digit approximations of the corresponding coordinates of the vertices, leaving out most of the details about how they were found and a more detailed proof that these correspond to exact embeddings. Those details are postponed to an upcoming longer paper. Nevertheless, with the help of a computer the data given in this note are sufficient to calculate arbitrary exact approximations of any of the presented embeddings from the defining algebraic equations, and thus to give convincing evidence that the Heawood graph indeed is unit distance embeddable, contrary to Chvatal's conjecture.  

One of these embeddings (the 9th in the list below) has been communicated to experts since May 2008 and has already been presented to a general mathematical audience at the Colloquium on Combinatorics in Magdeburg, Germany, in November 2008. The author only recently has learned of B.\ Horvat's PhD-thesis \cite{HorvatDiss}, dated April 2009, in which one further different from ours unit distance embedding of the Heawood graph was claimed to have been found. Since in that thesis the demonstration was presented in Slovenian only, we have not been able, yet, to compare results.

Special thanks go to Ed Pegg Jr., and Eric Weisstein (et al.\ from the technical support of Wolfram research), who in 2008 made available the computing facilities of Mathematica's at that time most recent iteration. Furthermore, for this research extensive use of the dynamical geometry software GeoGebra (\url{www.geogebra.org}) and the computer algebra system Singular (\url{www.singular.uni-kl.de}) was made. Fortunately both are still open domain, which should be honourably mentioned.
Finally the author gratefully acknowledges the generous financial support of Johannes Fuhs. Without any of the resources mentioned, this research would not have been possible.

\section{Determining an unit distance embedding for the Heawood graph, using Dynamical Geometry and numerics}

Let us first of all give some representation of the Heawood graph to fix notation: The Heawood graph is the point-line incidence graph of the smallest finite projective plane, the so-called Fano plane, i.e.\ the projective plane of order two, with seven points and seven lines. So if we let $P_i, 1\le i \le 7,$ denote the points and $l_i,1\le i\le 7,$ denote the lines of the plane, the 14 vertices of the Heawood graph will be given by the set $\{P_1,\dots,P_7,l_1,\dots,l_7\}.$ We specify the incidence between points and lines (and thus the adjacency relation of the graph) by  the following two pictures, which give the classical presentation of the Fano plane and its point-line incidence:  
\begin{figure}[h]
\begin{center}
$
\begin{array}{ccc}
{\includegraphics
[width=.42\linewidth]
{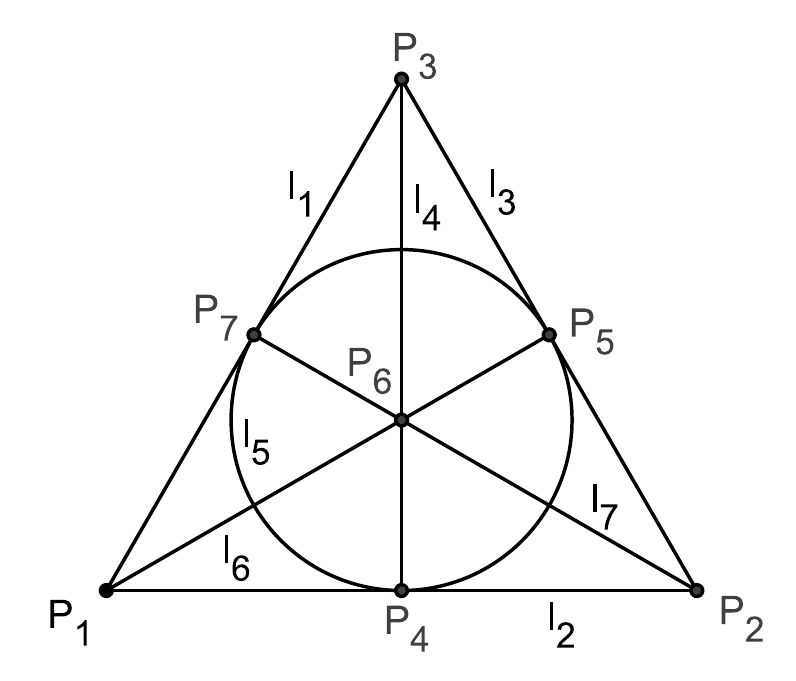}} 
&&
{\includegraphics[width=.42\linewidth]{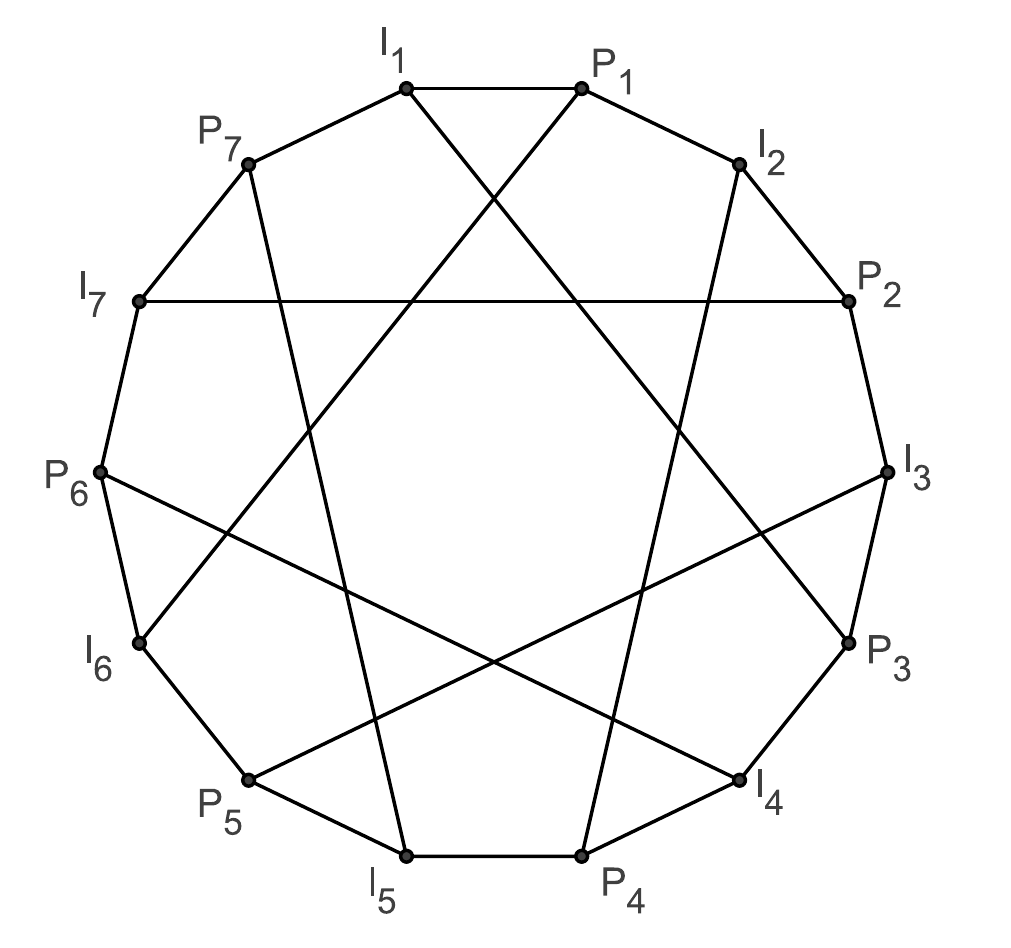}}\\
\hfill\hbox{(a)}\hfill &&\hfill\hbox{(b)}\hfill\, 
\end{array}
$
\caption{(a) The classic picture of the Fano plane, and (b) its incidence graph (i.e.\ the Heawood graph)}
\label{FanoHeawood}
\end{center}
\end{figure}

In order to determine a unit distance embedding of the Heawood graph, we need to find points $(x_{P_i},y_{P_i}),$ and $(x_{l_j},y_{l_j}),$ $1\le i,j\le 7,$ in $\mathbb{R}^2$ corresponding to the vertices of the graph such that $d(P_i,l_j) = 1$ holds if and only if $P_i$ is incident with $l_j.$
These constraints lead to 21 quadratic equations of the form
\begin{equation}
\label{DefEq}
(x_{P_i}-x_{l_j})^2+(y_{P_i}-y_{l_j})^2=1,
\end{equation}
one for each flag $(P_i\sim l_j)$ of the Fano plane. On the other hand we have $2\cdot 14 = 28$ (coordinate) variables, and thus a highly under-determined system of equations. Therefore we are allowed to fix an initial configuration of several vertices, and to proceed from there.
As suggested by Harris \cite{HarrisUDE}, we start by choosing a circle of maximal girth in the Heawood graph (which contains six elements). This we place in the shape of a rectangle, setting
\begin{equation}
\label{InitialConfig}
P_5:=(0,0),\, l_5:=(1,0),\, P_7=(1,1),\, l_7=(1,2),\, P_2=(0,2),\, l_3=(0,1). 
\end{equation}   

This leaves $28-12=16$ unknown variables, and $21-6=15$ equations. Thus we again may introduce one further restriction, which for the sake of computational simplicity, we choose to be that the points representing $l_4,$ $P_4,$ and $l_5$ should all lie on one straight line. This leads to $d(l_4,l_5)=2,$ which is the maximal possible distance between these two points. Thus the point corresponding to $l_4$ is supposed to lie on a circle of radius 2 around the origin of the drawing plane (which corresponds to $P_5$), and $P_4$ becomes the midpoint of the segment connecting $l_4$ and $l_5.$ 

With equations (\ref{DefEq}), and our special choices of coordinates in the initial configuration, see (\ref{InitialConfig}), plus the extra condition on $l_4\sim P_4\sim l_5,$ we have at hand a finite set of algebraic equations which completely determines a finite set of unit distance embeddings of the Heawood Graph. 

Let us explicitly list the set of equations complementing (\ref{InitialConfig}), using an order in which each point might be geometrically constructed from the previous ones by compass and ruler, after the position of e.g.\ $l_4$ has been fixed (the coordinates of any point $P$ are denoted by $(x_P,y_P)$):    
\begin{eqnarray}
(x_{l_4}-1)^2+y_{l_4}^2 - 4 & = & 0\label{onesegment}\\
x_{P_4}-\textstyle{\frac{1}{2}}(x_{l_4}+1) & = & 0\label{xP4}\\
y_{P_4}-\textstyle{\frac{1}{2}}y_{l_4} & = & 0\label{yP4}\\
x_{P_3}^2+(y_{P_3}-1)^2 - 1 & = & 0\label{P31}\\
(x_{P_3}-x_{l_4})^2+(y_{P_3}-y_{l_4})^2 - 1 & = & 0\label{P32}\\
(x_{P_6}-1)^2+(y_{P_6}-2)^2 - 1 & = & 0\label{P61}\\
(x_{P_6}-x_{l_4})^2+(y_{P_6}-y_{l_4})^2 - 1 & = & 0\label{P62}\\
x_{l_2}^2+(y_{l_2}-2)^2 - 1 & = & 0\label{l21}\\
(x_{l_2}-x_{P_4})^2+(y_{l_2}-y_{P_4})^2 - 1 & = & 0\label{l22}\\
(x_{l_1}-1)^2+(y_{l_1}-1)^2 - 1 & = & 0\label{l11}\\
(x_{l_1}-x_{P_3})^2+(y_{l_1}-y_{P_3})^2 - 1 & = & 0\label{l12}\\
x_{l_6}^2+y_{l_6}^2 - 1 & = & 0\label{l61}\\
(x_{l_6}-x_{P_6})^2+(y_{l_6}-y_{P_6})^2 - 1 & = & 0\label{l62}\\
(x_{P_1}-x_{l_2})^2+(y_{P_1}-y_{l_2})^2 - 1 & = & 0\label{P11}\\
(x_{P_1}-x_{l_6})^2+(y_{P_1}-y_{l_6})^2 - 1 & = & 0\label{P12}
\end{eqnarray}
\begin{eqnarray}%\\
(x_{P_1}-x_{l_1})^2+(y_{P_1}-x_{l_1})^2 - 1 & = & 0\label{unitdistanceCond}
\end{eqnarray}

Now, as will be shown in considerable detail in an upcoming (longer) paper, the set of equations (\ref{InitialConfig}) -- (\ref{unitdistanceCond}) completely determines a zero dimensional algebraic set of points in complex space $\mathbb{C}^{28},$ i.e., there are only finitely many solutions to these equations. 
In fact, the number of solutions is $79,$ only eleven of these are real and lead to eleven different unit distance embeddings of the Heawood graph (which all have the initial configuration and the fact that $l_4,$ $P_4,$ and $l_5$ lie on one line in common). The reason for there being only $79$ solutions is that equations (\ref{InitialConfig}) -- (\ref{unitdistanceCond}) lead to characteristic polynomials for each coordinate\footnote{The derivation of these was done in analogy to our approach in \cite{GerbrachtHarborth}, massively using the capabilities of Mathematica and Singular for calculating resultants and factoring polynomials. As an example, here we present the characteristic polynomial of the coordinate $x_{_{l_4}},$ which is 

\noindent
$\scriptstyle
p_{x_{l4}}(T)\,=\,3348011046054687446588586894387 + 273675328487397647237991825000783\, T +
    10528063279784456967398200502468691\, T^2 + 
    255652807673380729611728470237761555 \, T^3 +
    4422420653730204080254904433581059629 \, T^4 + 
    58239553681851019741523172701651095197 \, T^5 + 
    608930205226991194133708856923335926849 \, T^6 + 
    5203227805425306398124203036880713293545 \, T^7 + 
    37109973679879574898320679050599920287450 \, T^8 + 
    224472408717775611491021156521892619843158 \, T^9 + 
    1166012291532956694933924468283307736346382 \, T^{10} + 
    5253121604626527413065008160498494678879110 \, T^{11} + 
    20690863770430719393270631202371992513434414 \, T^{12} + 
    71715126275516155874072784490774971066237326 \, T^{13} + 
    219897164806211674807756610736580167553542758 \, T^{14 }+ 
    599083193386406195758633497190777431543886358 \, T^{15 }+ 
    1455265549140319863369871581645012065857955441 \, T^{16} + 
    3160933625571584072448347845721693351127774301 \, T^{17} + 
    6152912915312070842952691100801887803370907305 \, T^{18} + 
    10751995223766688842173817330681019107518783545 \, T^{19} + 
    16888459659695355326863471817880692818622047623 \, T^{20} + 
    23863989284324858347511498529857889181323950967 \, T^{21} + 
    30346538554876120431728853077314517314609386819 \, T^{22} + 
    34722813139066795200081139797717329223025992699 \, T^{23} + 
    35716781564909427260214641236767872783162088204 \, T^{24} + 
    32963773017875955980864755706102737727974961688 \, T^{25} + 
    27201188778043412156622512508710379868716241416 \, T^{26} + 
    19954407479150801176386566760213350973570196080 \, T^{27} + 
    12902691890291653798206974719870993995735753540 \, T^{28} + 
    7274584518541872070335933586256322019748139172 \, T^{29} + 
    3550298683130683434462662943215234037891507412 \, T^{30} + 
    1533983381070251025995839971747580678500964852 \, T^{31} + 
    664103288660372783854070699409333594554864741 \, T^{32} + 
    355269696471069385886754716351566237266830009 \, T^{33} + 
    213238754173051016042819729417269617854966165 \, T^{34} + 
    77130998985650655864689962089382720577858101 \, T^{35 }- 
    57662664820854923809493690824194000968797973 \, T^{36 }- 
    146045321267662575006252144965793560225509061 \, T^{37} - 
    164022007275895644197052849670790737036540873 \, T^{38} - 
    126213399593210124294769323126921742027688497 \, T^{39} - 
    67869160289243415287139367956058055810404822 \, T^{40} - 
    19570606574427556470966233073766236628787234 \, T^{41} + 
    6140751881298455069763046326781936407849238 \, T^{42 }+ 
    12720991312674958659494390034544200285598942 \, T^{43} + 
    9840137643451726574992603743314811193317886 \, T^{44} + 
    5180867575272248126071836905848828341927070 \, T^{45} + 
    1923952833473734147443634652898764867278198 \, T^{46} + 
    286935408276107233753158822122577885606822 \, T^{47} - 
    343872926425618669220741202688368202345065 \, T^{48} - 
    451645674349824891937650532097542435080453 \, T^{49} - 
    325157218431048323421805399697113970403121 \, T^{50} - 
    152272756904971138353148344210050406803617 \, T^{51} - 
    30934416501269415569285918492882277029311 \, T^{52} + 
    21867253654523569285667250014704999794577 \, T^{53} + 
    28955348159492426037443729536713509636773 \, T^{54} + 
    17321709733106215547946139735151891780269 \, T^{55} + 
    4733784662326174469816987234959768253776 \, T^{56} - 
    1650827959998751884275421145646879272940 \, T^{57} - 
    2435243231716218466580115477132980137292 \, T^{58} - 
    1097279690260575519531876572540059803892 \, T^{59} - 
    60617631026953339799378305296984824656 \, T^{60} + 
    200727376265061817580032667984094835280 \, T^{61} + 
    109385892925207478360122518287948266224 \, T^{62} + 
    14201705397119143149709337683063717104 \, T^{63} - 
    11463391775661584618715895715025904128 \, T^{64} - 
    6556557400356413683063078157405200320 \, T^{65} - 
    898635822877066299154282762314323520 \, T^{66} + 
    477056183905245971917488031692938304 \, T^{67} + 
    254616663098419271111012531383618560 \, T^{68} + 
    31343179682405215504161837658819584 \, T^{69} - 
    14303114368662112977785429692643328 \, T^{70} - 
    6892489761595983453459595854256128 \, T^{71} - 
    763800345871643605733535512788992 \, T^{72} + 
    341989984727973884867396338188288 \, T^{73 }+ 
    149189048927171391219263917572096 \, T^{74} + 
    11025799477301561380923592949760 \, T^{75} - 
    8175821639408563679884718899200 \, T^{76} - 
    2120259444356145889512456192000 \, T^{77} + 
    152135800369825007098920960000 \, T^{78} + 
    82521703002365615643033600000 \, T^{79}.
$}
 of degree $79,$ which each have eleven real zeroes. These can be grouped, again, to give the eleven different real solutions.

\eject
\begin{theorem}
Numerical approximations (up to 15 digits) to the real solutions of the above set of equations with corresponding embeddings are  given by

\smallskip
%\eject
\begin{minipage}{.5\linewidth}
\hspace{-1truecm}
{\includegraphics[width=\linewidth]{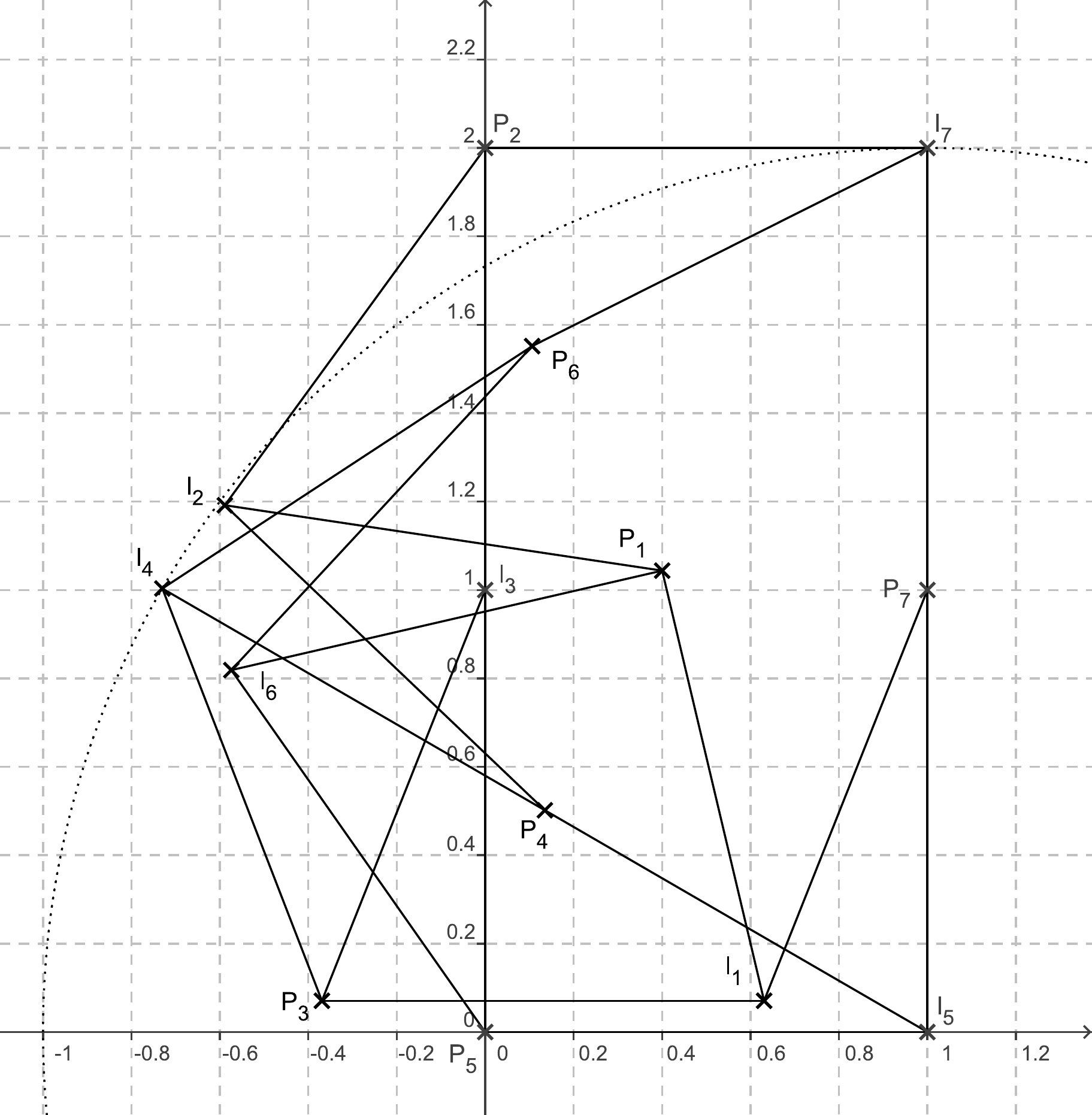}}
\hspace{.5truecm} 
\end{minipage}
%&&
\begin{minipage}{.45\linewidth}
$
\vspace{-4truecm}
\begin{array}{lcl}%% A6 %%
P_1 &\simeq & (0.400182002388641,\, 1.043762692468723)\nonumber\\
P_3 &\simeq & (-0.369307668700666,\, 0.070692814060453)\nonumber\\
P_4 &\simeq & (0.134937917545110,\, 0.501664821866961)\nonumber\\
P_6 &\simeq & (0.106134457655163,\, 1.551664866189844)\nonumber\\
l_1 &\simeq & (0.630692331299334,\, 0.070692814060453)\nonumber\\
l_2 &\simeq & (-0.588810425254328,\, 1.191728830706045)\nonumber\\
l_4 &\simeq & (-0.730124164909779,\, 1.003329643733922)\nonumber\\
l_6 &\simeq & (-0.574170534719569,\, 0.818735730904572)\nonumber;\\
\vspace{5truecm}
\end{array}
%\end{array}
$
\end{minipage}

\smallskip
\begin{minipage}{.5\linewidth}
\hspace{-1truecm}
{\includegraphics[width=\linewidth]{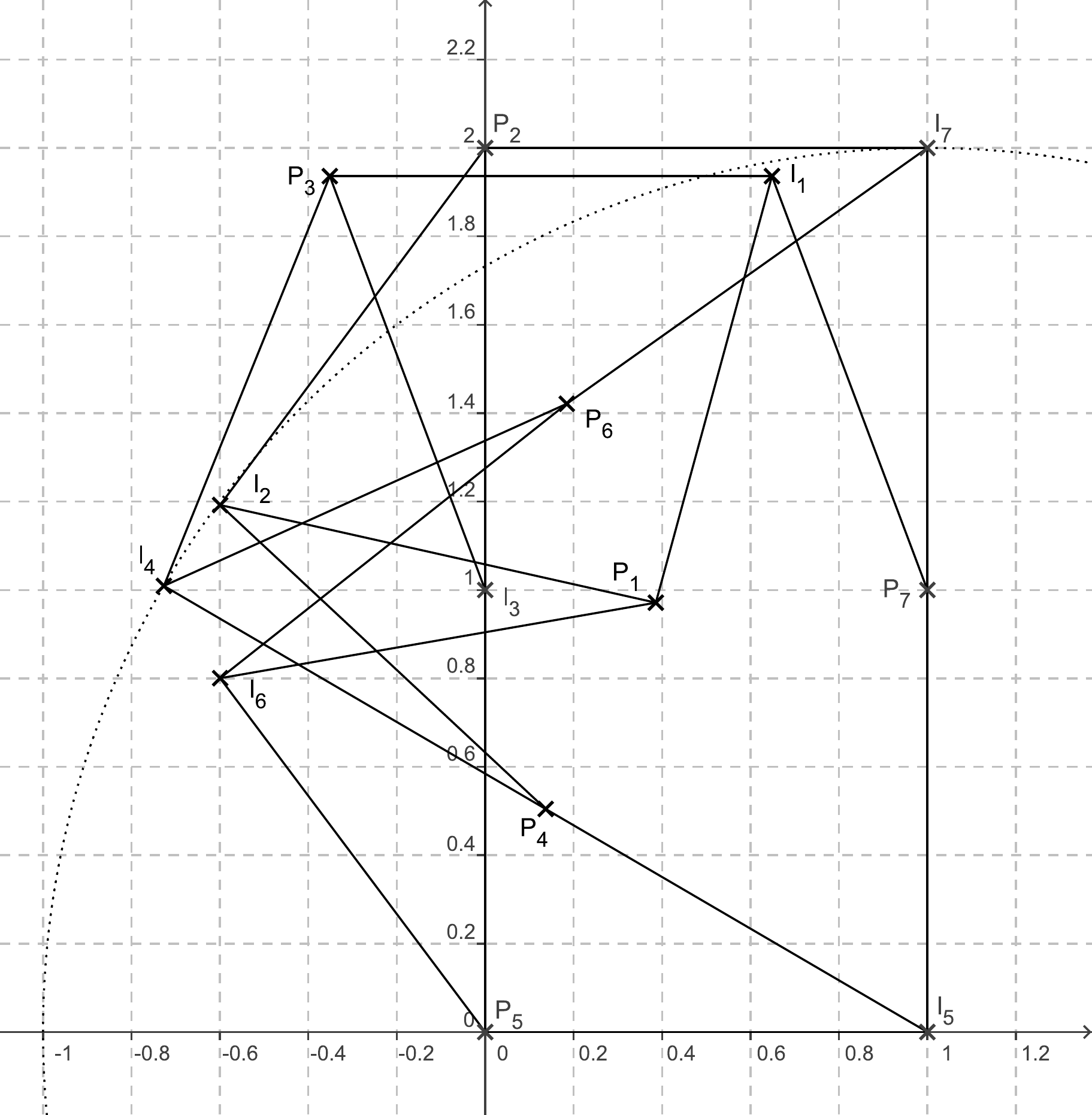}}
\hspace{.5truecm} 
\end{minipage}
%&&
\begin{minipage}{.45\linewidth}
$
\vspace{-4truecm}
\begin{array}{lcl}%% A7 %%
P_1 &\simeq & (0.385844376323838,\, 0.971300460792625)\nonumber\\
P_3 &\simeq & (-0.351496569091080,\, 1.936189169942272)\nonumber\\
P_4 &\simeq & (0.136658400253077,\, 0.504619938316376)\nonumber\\
P_6 &\simeq & (0.184497979349687,\, 1.421245773825144)\nonumber\\
l_1 &\simeq & (0.648503430908920,\, 1.936189169942272)\nonumber\\
l_2 &\simeq & (-0.589452825277335,\, 1.192197198090668)\nonumber\\
l_4 &\simeq & (-0.726683199493847,\, 1.009239876632752)\nonumber\\
l_6 &\simeq & (-0.599446494990155,\, 0.800414829725199)\nonumber;\\
\vspace{5truecm}
\end{array}
%\end{array}
$
\end{minipage}

\smallskip
\begin{minipage}{.5\linewidth}
\hspace{-1truecm}
{\includegraphics[width=\linewidth]{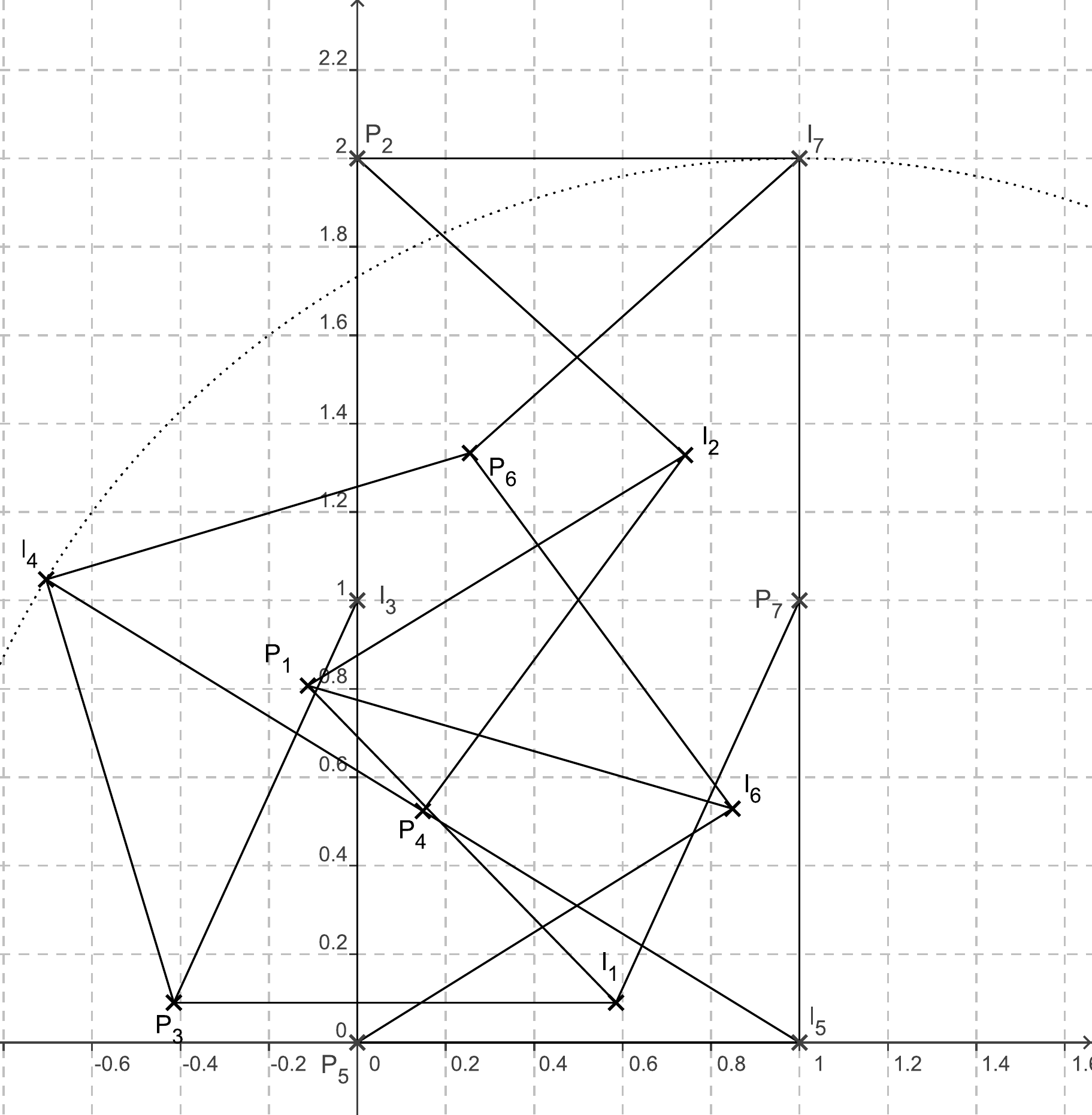}}
\hspace{.5truecm}
\end{minipage} 
%&&
\begin{minipage}{.45\linewidth}
$
%\vspace{-3truecm}
\begin{array}{lcl}%% A5 %%
P_1 &\simeq & (-0.111888421172288,\, 0.807281254230185)\nonumber\\
P_3 &\simeq & (-0.414946144051627,\, 0.090154025377544)\nonumber\\
P_4 &\simeq & (0.148144628769197,\, 0.523777077109366)\nonumber\\
P_6 &\simeq & (0.254555103623971,\, 1.333432744228363)\nonumber\\
l_1 &\simeq & (0.585053855948373,\, 0.090154025377544)\nonumber\\
l_2 &\simeq & (0.741321136264131,\, 1.328849515437814)\nonumber\\
l_4 &\simeq & (-0.703710742461606,\, 1.047554154218733)\nonumber\\
l_6 &\simeq & (0.848614578463299,\, 0.529011622953180)\nonumber;\\
%\vspace{4truecm}
\end{array}
%\end{array}
$
\end{minipage}

\smallskip
\begin{minipage}{.5\linewidth}
\hspace{-1truecm}
{\includegraphics[width=\linewidth]{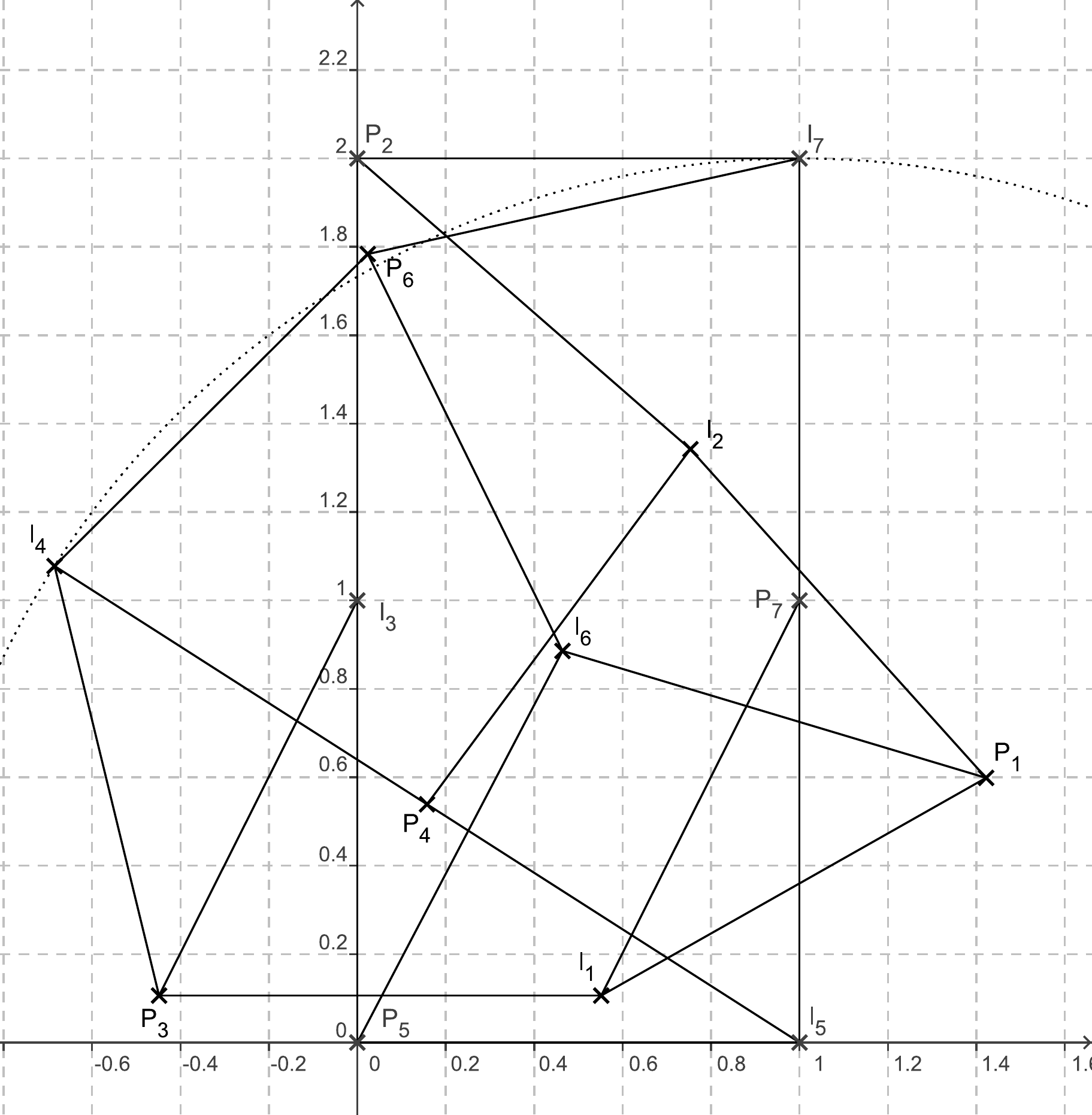}}
\hspace{.5truecm}
\end{minipage} 
%&&
\begin{minipage}{.45\linewidth}
$
%\vspace{-3truecm}
\begin{array}{lcl}%% A4 %%
P_1 &\simeq & (1.421904368600474,\, 0.598683521299139)\nonumber\\
P_3 &\simeq & (-0.448398217164224,\, 0.106166101088158)\nonumber\\
P_4 &\simeq & (0.157643970813835,\, 0.538921441486342)\nonumber\\
P_6 &\simeq & (0.023681878654164,\, 1.783660161015737)\nonumber\\
l_1 &\simeq & (0.551601782835776,\, 0.106166101088158)\nonumber\\
l_2 &\simeq & (0.753214201921679,\, 1.342224684239756)\nonumber\\
l_4 &\simeq & (-0.684712058372330,\, 1.077842882972684)\nonumber\\
l_6 &\simeq & (0.464016631427877,\, 0.885826487388092)\nonumber;\\
%\vspace{4truecm}
\end{array}
%\end{array}
$
\end{minipage}

\smallskip
\begin{minipage}{.5\linewidth}
\hspace{-1truecm}
{\includegraphics[width=\linewidth]{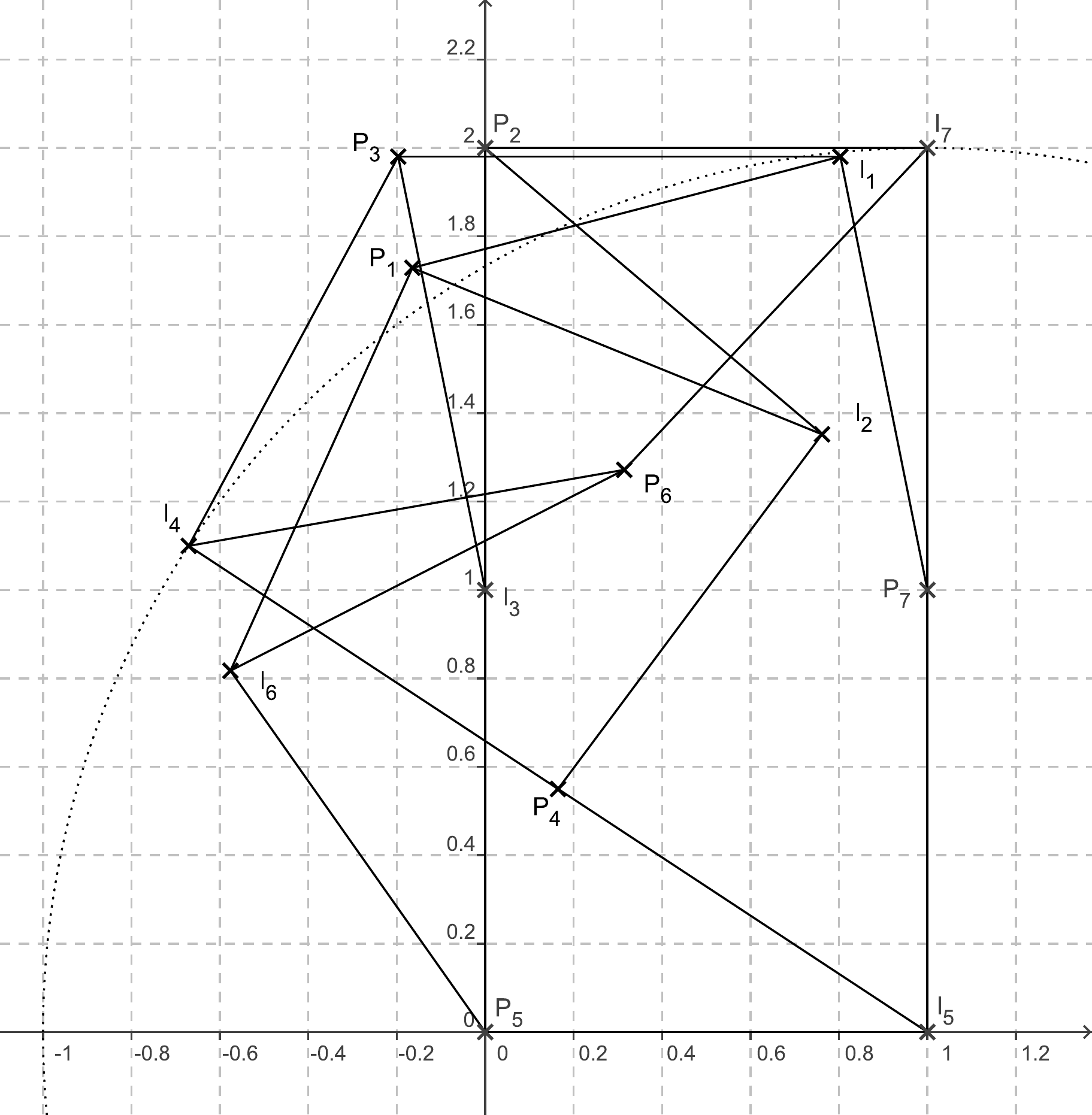}}
\hspace{.5truecm}
\end{minipage} 
%&&
\begin{minipage}{.45\linewidth}
$
%\vspace{-3truecm}
\begin{array}{lcl}%% A8 %%
P_1 &\simeq & (-0.164629969634109,\, 1.728731712593544)\nonumber\\
P_3 &\simeq & (-0.196826391803828,\, 1.980438356802449)\nonumber\\
P_4 &\simeq & (0.164749301040245,\, 0.549869320736518)\nonumber\\
P_6 &\simeq & (0.314574903717282,\, 1.271856856527629)\nonumber\\
l_1 &\simeq & (0.803173608196172,\, 1.980438356802449)\nonumber\\
l_2 &\simeq & (0.761740164689577,\, 1.352117355149333)\nonumber\\
l_4 &\simeq & (-0.670501397919511,\, 1.099738641473037)\nonumber\\
l_6 &\simeq & (-0.576161224535975,\, 0.817336065117162)\nonumber;\\
%\vspace{4truecm}
\end{array}
%\end{array}
$
\end{minipage}

\smallskip
\begin{minipage}{.5\linewidth}
\hspace{-1truecm}
{\includegraphics[width=\linewidth]{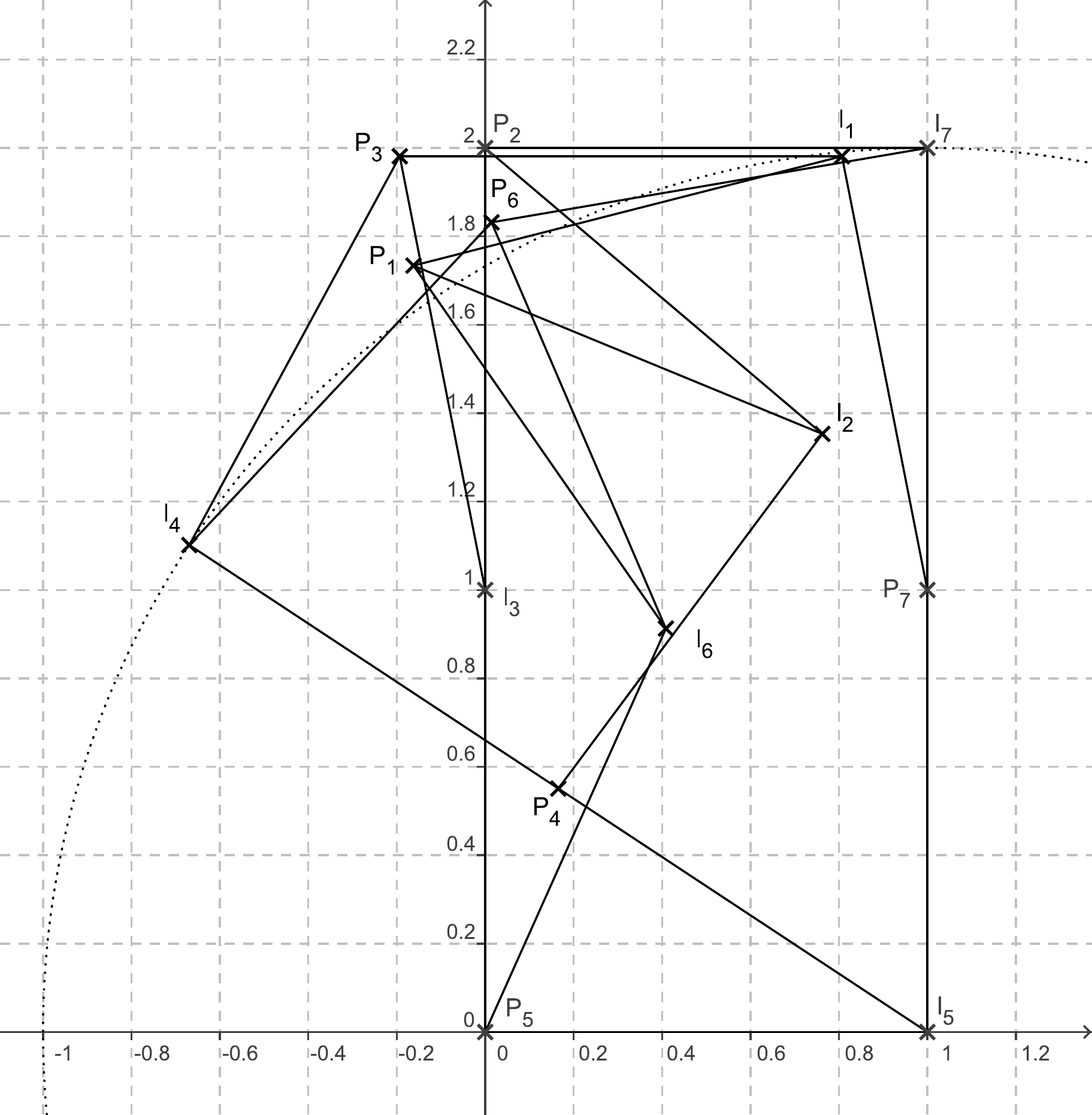}}
\hspace{.5truecm}
\end{minipage} 
%&&
\begin{minipage}{.45\linewidth}
$
%\vspace{-3truecm}
\begin{array}{lcl}%% A9 %%
P_1 &\simeq & (-0.162159171174261,\, 1.733808025238098)\nonumber\\
P_3 &\simeq & (-0.193231727879694,\, 1.981153147750456)\nonumber\\
P_4 &\simeq & (0.165394595983793,\, 0.550848272745721)\nonumber\\
P_6 &\simeq & (0.014270178593230,\, 1.831664860503289)\nonumber\\
l_1 &\simeq & (0.806768272120306,\, 1.981153147750456)\nonumber\\
l_2 &\simeq & (0.762499622579621,\, 1.353011340465742)\nonumber\\
l_4 &\simeq & (-0.669210808032414,\, 1.101696545491441)\nonumber\\
l_6 &\simeq & (0.408620165537735,\, 0.912704530675680)\nonumber;\\
%\vspace{4truecm}
\end{array}
%\end{array}
$
\end{minipage}

\smallskip
\begin{minipage}{.5\linewidth}
\hspace{-1truecm}
{\includegraphics[width=\linewidth]{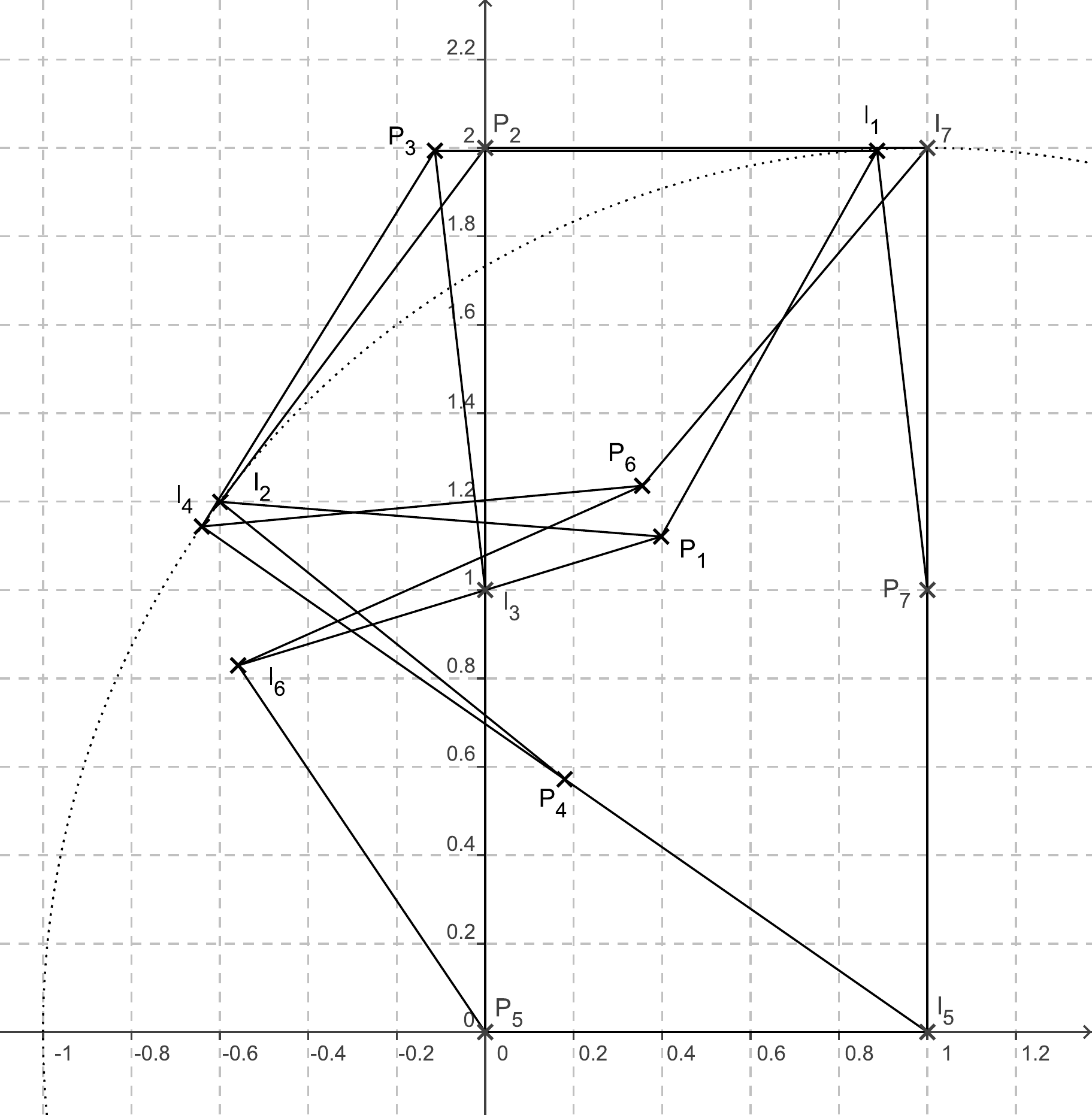}}
\hspace{.5truecm}
\end{minipage} 
%&&
\begin{minipage}{.45\linewidth}
$
%\vspace{-3truecm}
\begin{array}{lcl}%% A10 %%
P_1 &\simeq & (0.39788958050138643,\, 1.120875683482126)\nonumber\\
P_3 &\simeq & (-0.11373768509513692,\, 1.993510814731878)\nonumber\\
P_4 &\simeq & (0.17962533307815014,\, 0.571826377384660)\nonumber\\
P_6 &\simeq & (0.35499397395823023,\, 1.235822516446733)\nonumber\\
l_1 &\simeq & (0.88626231490486308,\, 1.993510814731878)\nonumber\\
l_2 &\simeq & (-0.59903243717863292,\, 1.199275241292098)\nonumber\\
l_4 &\simeq & (-0.64074933384369972,\, 1.143652754769321)\nonumber\\
l_6 &\simeq & (-0.55868296106963737,\, 0.829381304955967)\nonumber;\\
%\vspace{4truecm}
\end{array}
%\end{array}
$
\end{minipage}

\smallskip
\begin{minipage}{.5\linewidth}
\hspace{-1truecm}
{\includegraphics[width=\linewidth]{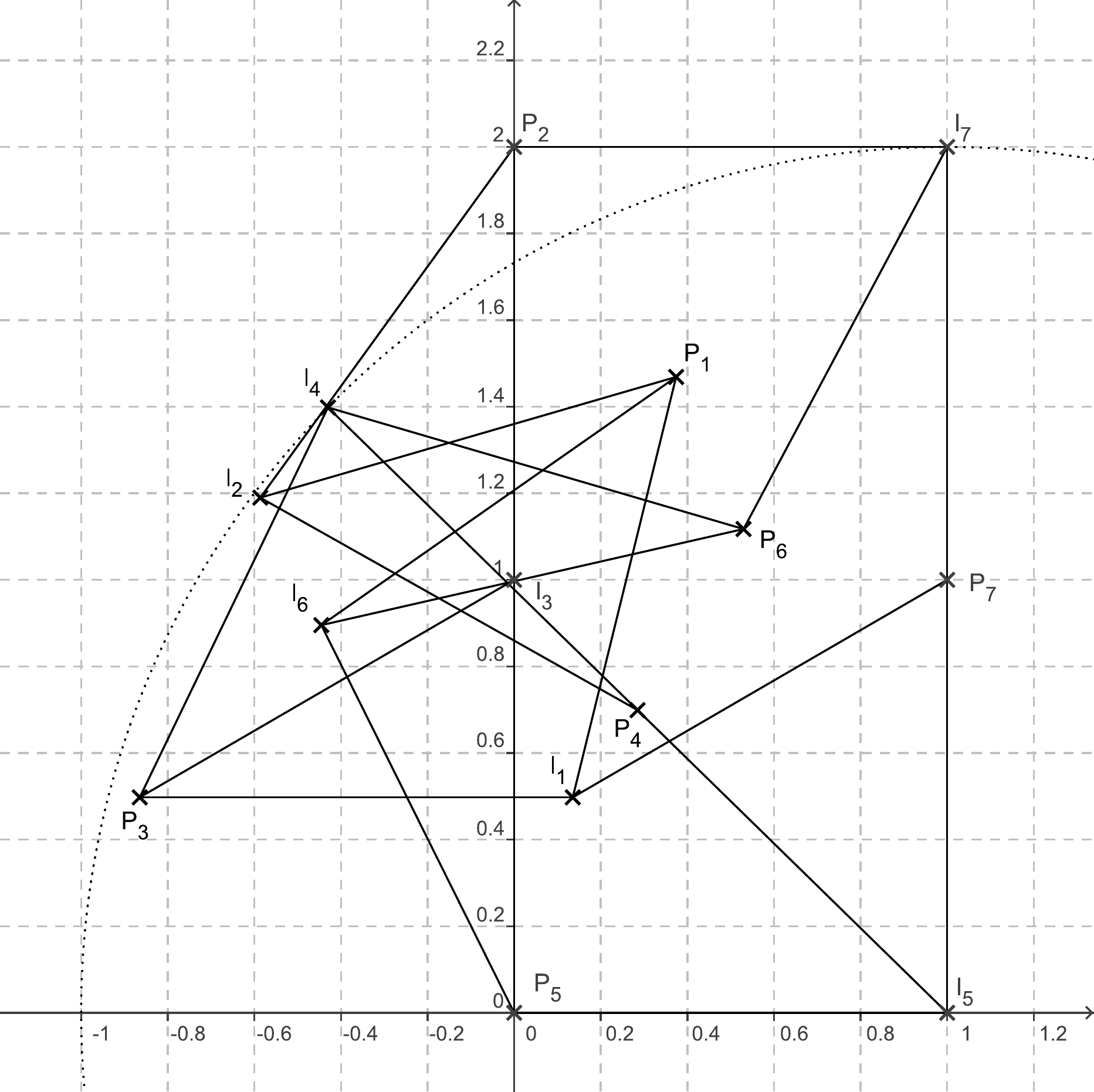}}
\hspace{.5truecm}
\end{minipage} 
%&&
\begin{minipage}{.45\linewidth}
$
%\vspace{-3truecm}
\begin{array}{lcl}%% A3 %%
P_1 &\simeq & (0.374047478003772,\, 1.468744576326795)\nonumber\\
P_3 &\simeq & (-0.864667172852078,\, 0.497654819678744)\nonumber\\
P_4 &\simeq & (0.284999030498420,\, 0.699123460922175)\nonumber\\
P_6 &\simeq & (0.529767447111884,\, 1.117457453601060)\nonumber\\
l_1 &\simeq & (0.135332827147922,\, 0.497654819678744)\nonumber\\
l_2 &\simeq & (-0.586287978493963,\, 1.189897286590483)\nonumber\\
l_4 &\simeq & (-0.430001939003160,\, 1.398246921844351)\nonumber\\
l_6 &\simeq & (-0.445265782381637,\, 0.895398449317436)\nonumber;\\
%\vspace{4truecm}
\end{array}
%\end{array}
$
\end{minipage}

\smallskip
\begin{minipage}{.5\linewidth}
\hspace{-1truecm}
{\includegraphics[width=\linewidth]{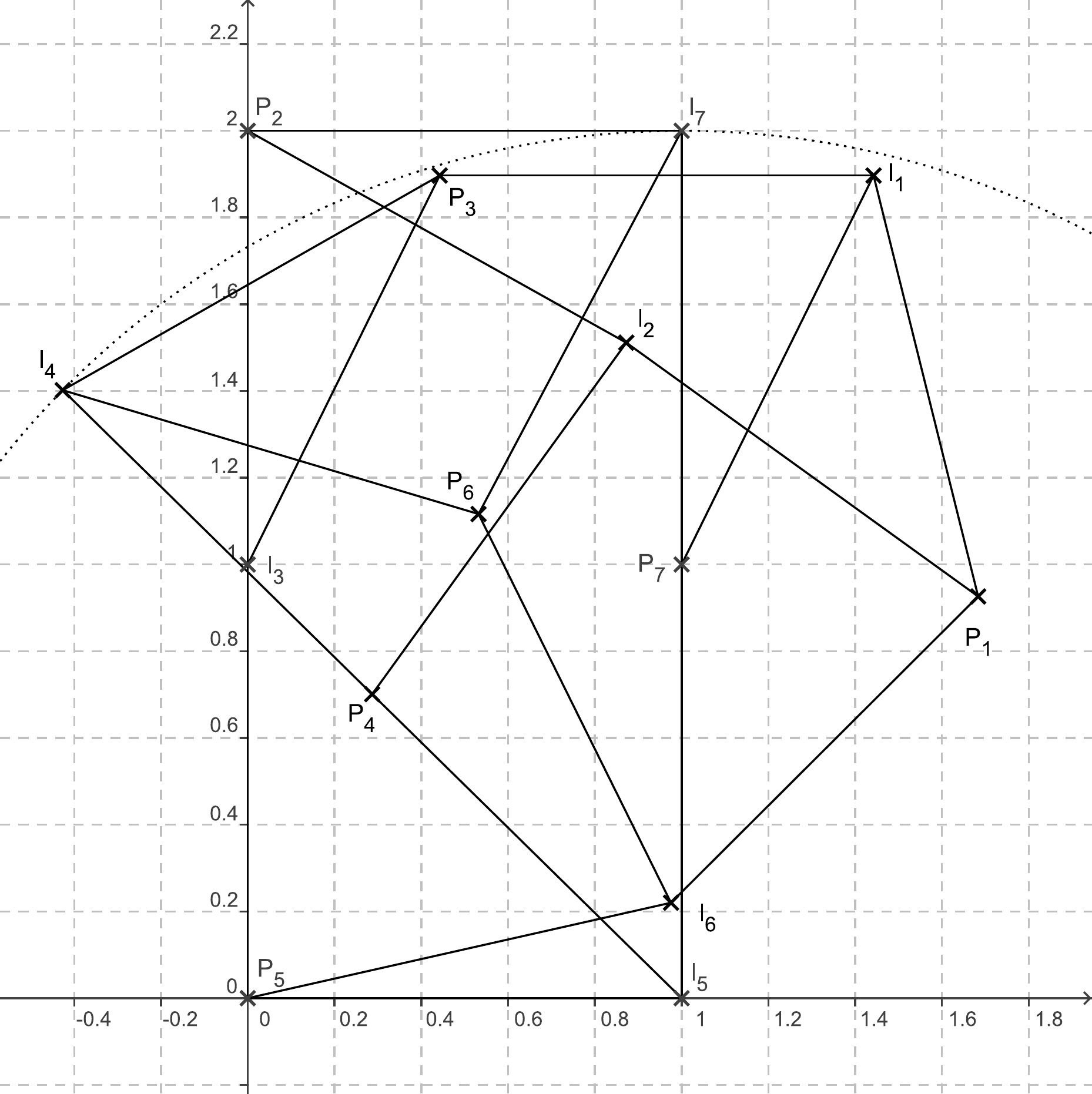}}
\hspace{.5truecm}
\end{minipage} 
%&&
\begin{minipage}{.45\linewidth}
$
%\vspace{-3truecm}
\begin{array}{lcl}%% A11 %%
P_1 &\simeq & (1.683478977241937,\, 0.926313519981179)\nonumber\\
P_3 &\simeq & (0.442398409684554,\, 1.896818625599149)\nonumber\\
P_4 &\simeq & (0.286751573157956,\, 0.700911322217696)\nonumber\\
P_6 &\simeq & (0.531884805746619,\, 1.116332548454523)\nonumber\\
l_1 &\simeq & (1.442398409684554,\, 1.896818625599149)\nonumber\\
l_2 &\simeq & (0.872507318643794,\, 1.511398957306269)\nonumber\\
l_4 &\simeq & (-0.426496853684088,\, 1.401822644435391)\nonumber\\
l_6 &\simeq & (0.975476510630742,\, 0.220103560214309)\nonumber;\\
%\vspace{4truecm}
\end{array}
%\end{array}
$
\end{minipage}

\smallskip
\begin{minipage}{.5\linewidth}
\hspace{-1truecm}
{\includegraphics[width=\linewidth]{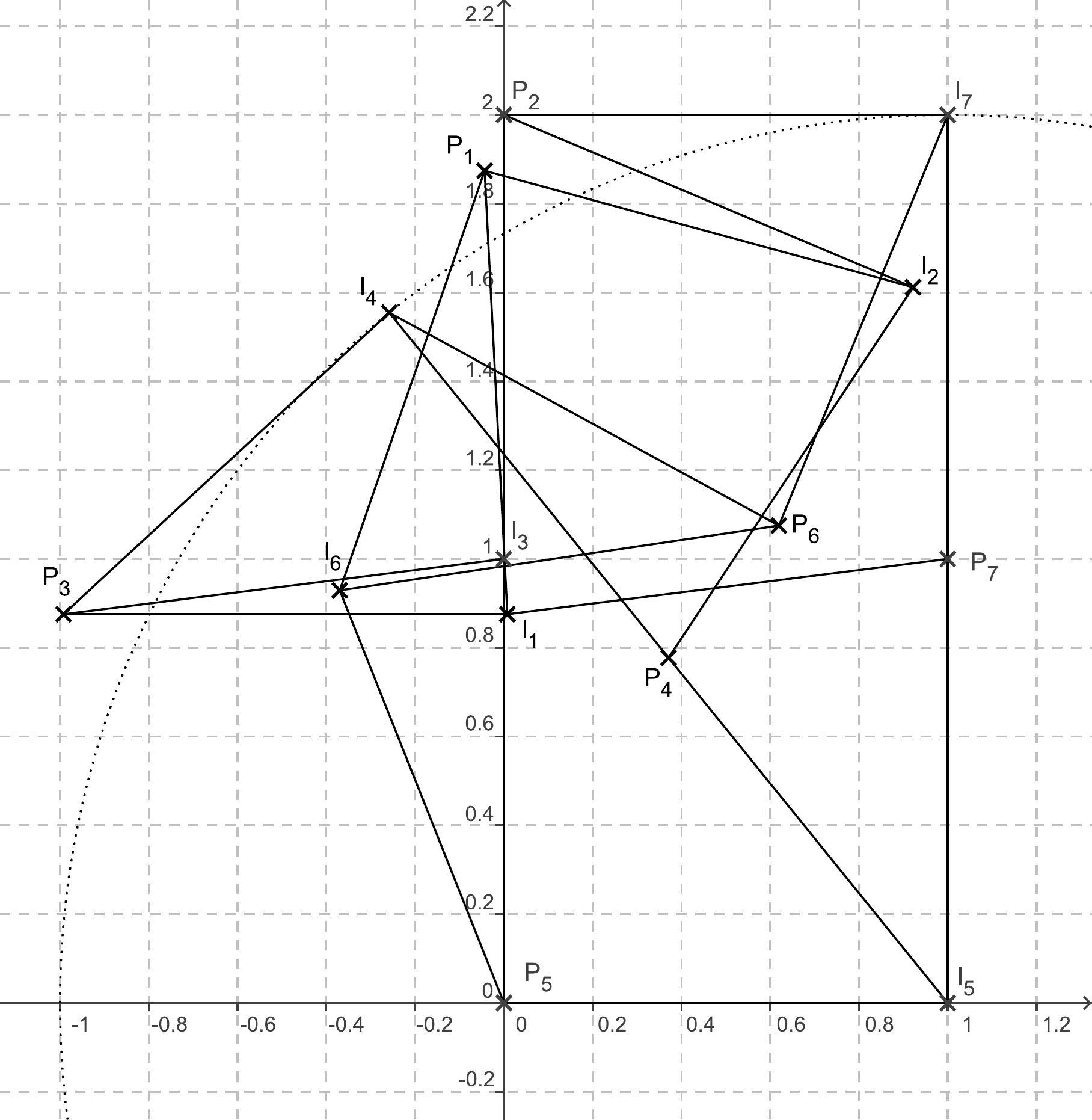}}
\hspace{.5truecm}
\end{minipage} 
%&&
\begin{minipage}{.45\linewidth}
$
%\vspace{-3truecm}
\begin{array}{lcl}%% A1 %%
P_1 &\simeq & (-0.043329551699449,\, 1.874285721361077)\nonumber\\
P_3 &\simeq & (-0.992231159457987,\, 0.875592097515236)\nonumber\\
P_4 &\simeq & (0.370915641186566,\, 0.777337037260087)\nonumber\\
P_6 &\simeq & (0.619416519241032,\, 1.075253432461959)\nonumber\\
l_1 &\simeq & (0.007768840542013,\, 0.875592097515236)\nonumber\\
l_2 &\simeq & (0.921663138782111,\, 1.612008945192925)\nonumber\\
l_4 &\simeq & (-0.258168717626869,\, 1.554674074520175)\nonumber\\
l_6 &\simeq & (-0.369844480601252,\, 0.929093676745671)\nonumber;\\
%\vspace{4truecm}
\end{array}
%\end{array}
$
\end{minipage}

\smallskip
\begin{minipage}{.5\linewidth}
\hspace{-1truecm}
{\includegraphics[width=\linewidth]{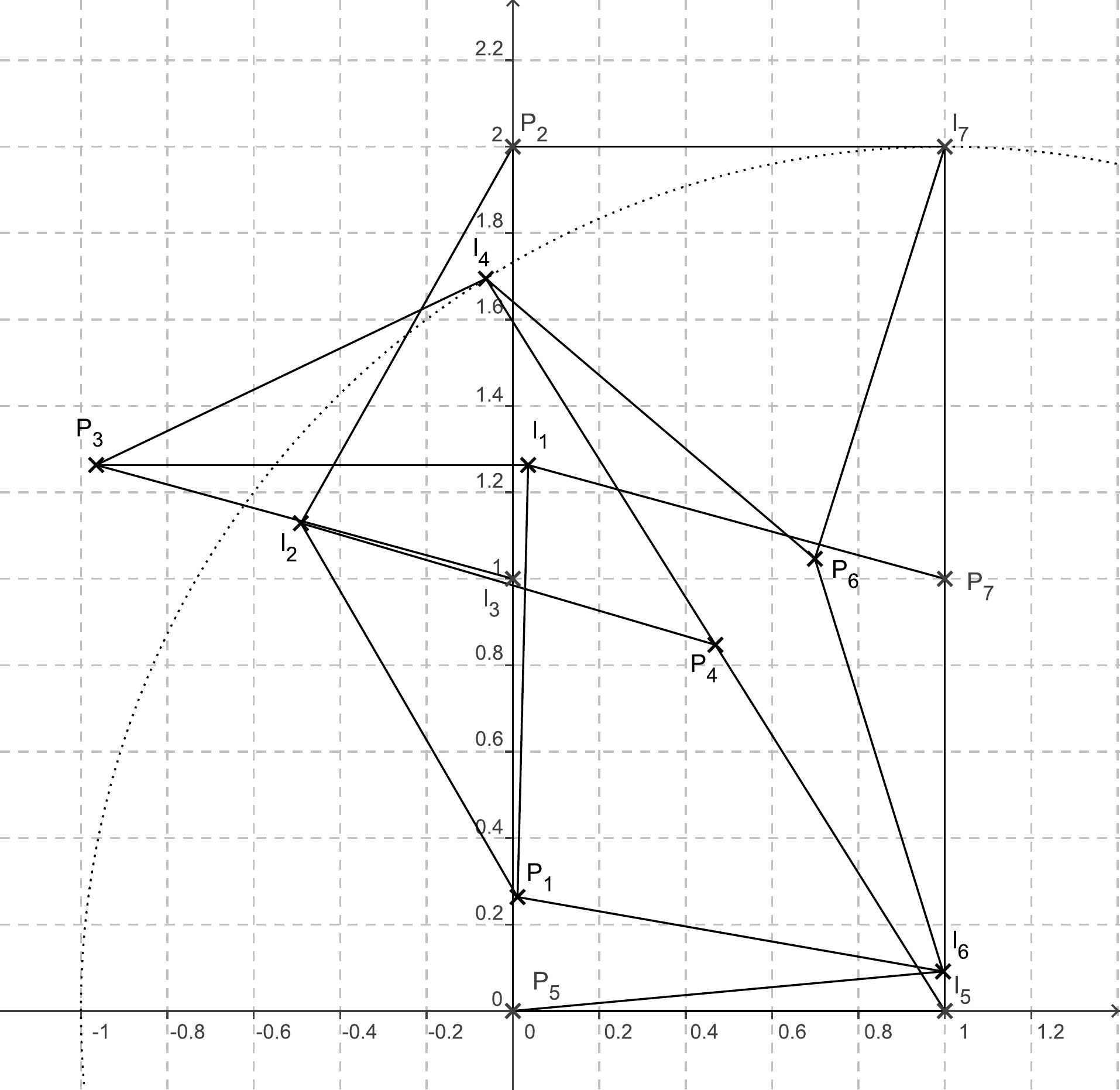}}
\hspace{.5truecm}
\end{minipage} 
%&&
\begin{minipage}{.45\linewidth}
$
%\vspace{-3truecm}
\begin{array}{lcl}%% A2 %%
P_1 &\simeq & (0.010754855391715,\, 0.263588750091565)\nonumber\\
P_3 &\simeq & (-0.964717307331239,\, 1.263287897434660)\nonumber\\
P_4 &\simeq & (0.468775634039814,\, 0.847231180381246)\nonumber\\
P_6 &\simeq & (0.699093018450262,\, 1.046346505037272)\nonumber\\
l_1 &\simeq & (0.035282692668761,\, 1.263287897434660)\nonumber\\
l_2 &\simeq & (-0.490788504254845,\, 1.128721259245187)\nonumber\\
l_4 &\simeq & (-0.062448731920371,\, 1.694462360762491)\nonumber\\
l_6 &\simeq & (0.995815833199486,\, 0.091382855882344)\nonumber.
%\vspace{4truecm}
\end{array}
%\end{array}
$
\end{minipage}
\end{theorem}

As a final remark we would like to state the observation that all of these embeddings are {\sl regular,} i.e., embedded vertices only coincide with the embeddings of those edges, which they are incident with. 

%------------------------------------------------------------------------------

\bibliographystyle{amsplain}
\bibliography{Heawood}
\end{document}